\documentclass[12pt]{amsart}

\usepackage{graphicx}
\usepackage{amsthm}
\usepackage{amsmath, amssymb}
\swapnumbers

\theoremstyle{plain}
\newtheorem{Thm}{Theorem}

\newtheorem{Coro}[Thm]{Corollary}
\newtheorem{Lem}[Thm]{Lemma}

\theoremstyle{definition}
\newtheorem{Def}[Thm]{Definition}

\begin{document}

\title{Layered models for closed 3-manifolds}

\author{Jesse Johnson}
\address{\hskip-\parindent
        Department of Mathematics \\
        Oklahoma State University \\
        Stillwater, OK 74078 \\
        USA}
\email{jjohnson@math.okstate.edu}


%\thanks{Research supported by NSF MSPRF grant 0602368}

\begin{abstract}
We define a combinatorial structure on 3-manifolds that combines the model manifolds constructed in Minsky's proof of the ending lamination conjecture with the layered triangulations defined by Jaco and Rubinstein.
\end{abstract}

\maketitle

Minsky's recent proof of the ending lamination conjecture~\cite{minsky} introduced a construction in which a path in the pants complex of a boundary component $S$ of a hyperbolic 3-manifold $M$ is used to build a combinatorial model of a hyperbolic structure on $S \times [0,\infty)$.  This structure is analogous in many ways to a triangulation:  It consists of a one-dimensional object (a link) spanned by two-dimensional pieces (embedded pairs of pants) that cut the manifold into relatively simple pieces (which we will describe below).

The construction of these models is closely related to the layering construction used to build ideal triangulations of punctured torus bundles in~\cite{guer}.  In particular, a model is determined by a path hierarchy in the curve complex, or (equivalently) a path in the pants complex for the boundary.  The layered triangulations in~\cite{guer} correspond to path hierarchies in the arc complex for the punctured surface.

Jaco and Rubinstein~\cite{jacoru} have studied layered triangulations of handlebodies, which can be glued together to form layered triangulations of closed 3-manifolds.  Applying the analogy between Minsky models and layered triangulations suggests a way to construct Minsky models for closed 3-manifolds.  In the present paper, we study this construction, applying results from the theory of Heegaard splittings to prove the following:

\begin{Thm}
\label{mainthm}
Every compact, hyperbolic 3-manifold admits a model decomposition in which every component of the 1-skeleton is a knotted loop.  % any the complement of the link is hyperbolic
\end{Thm}

By a \textit{knotted} loop in an arbitrary 3-manifold, we mean a loop that does not bound an embedded disk.  Such a loop may still be homotopy trivial in $M$.  

This construction can be readily extended to 3-manifolds with atoroidal boundary.  It should also be possible to extend it to 3-manifolds with toroidal boundary by modifying the definition of model blocks slightly to allow them to ``wrap'' around the torus boundary, which is more in line with the construction in~\cite{minsky}.  However, to simplify the exposition, we will not consider these generalizations here.

We also note the following property of any model of a closed 3-manifold.  This result does not require that the 3-manifold be hyperbolic, or that the model structure satisfy the conclusions of Theorem~\ref{mainthm}.

\begin{Thm}
\label{mainthm2}
The 1-skeleton of any model structure for a compact, closed, orientable 3-manifold is a hyperbolic link.
\end{Thm}

The outline of the paper is as follows:  We define the model structure in Section~\ref{modelsect}, then review the theory of generalized Heegaard splittings in Section~\ref{genhsect}.  We define a general construction of models for compression bodies (the building blocks of generalized Heegaard splittings) in Section~\ref{compsect}, then in Section~\ref{modelcompsect} we refine this construction to produce models for compression bodies with no unknotted loops.  Finally, we use this refined construction to prove Theorem~\ref{mainthm} in Section~\ref{constrsect}.

\section{Model structures}
\label{modelsect}

A \textit{pants decomposition} for a compact, closed, orientable surface $S$ is a set $P$ of pairwise disjoint, essential simple closed curves in $S$ such that each component of $S \setminus P$ is a thrice-punctured sphere (i.e. a pair of pants).

Given a pants decomposition $P$ and a loop $\ell \subset P$, the complement $S \setminus (P \setminus \ell)$ consists of a collection of pairs of pants and either a four-punctured sphere or a once-punctured torus $F_\ell \subset S$ that contains $\ell$.  When $F_\ell$ is a punctured torus and $\ell'$ is a second essential loop in $F_\ell$ not parallel to $\ell$, we will say that $\ell'$ intersects $\ell$ minimally if  $\ell \cap \ell'$ is a single point.  When $F_\ell$ is a four-punctured sphere and $\ell' \subset F_\ell$ is an essential, non-parallel loop we will say $\ell'$ intersects $\ell$ minimally if $\ell \cap \ell'$ is exactly two points.  In both cases, replacing $\ell$ with a new loop that intersects $\ell$ minimally produces a new pants decomposition for $S$.  (The condition that $\ell'$ intersects $\ell$ minimally is not necessary to produce a new pants decompositon, but it is useful for other reasons.)

\begin{Def}
The \textit{pants complex} $\mathcal{P}(S)$ for a compact, closed surface $S$ is the graph whose vertices are isotopy classes of pants decompositions for $S$ and whose edges connect each pants decomposition $P$ to each pants decomposition that results from replacing a loop of $P$ with a new essential loop in $F_\ell$ that intersects $\ell$ minimally.
\end{Def}

An edge path in the pants complex for $S$ determines a sequence of pants decompositions for $S$.  We will use this sequence of pants decompositions to construct a combinatorial structure on $S \times [0,1]$ as follows:

\begin{Def}
A \textit{pants block} is a handlebody with a collection of essential loops in its boundary, forming a pants decomposition of one of two forms:  The first type can be identified with a once punctured torus cross an interval $T \times [0,1]$ so that one loop is essential in the annulus $\partial T \times [0,1]$ and the other two loops are contained in $T \times \{0\}$ and $T \times \{1\}$ so that their projections into $T$ intersect in a single point.  The second type can be identified with a four-times punctured sphere cross an interval $T \times [0,1]$ so that one loop is essential in each component of $\partial T \times [0,1]$ and the remaining two loops are contained in $T \times \{0\}$ and $T \times \{1\}$ so that their projections into $T$ intersect in exactly two points.  
\end{Def}

In the case when $T$ is a once punctured torus, the complement of the loops is two pairs of pants.  For $T$ a four punctured sphere, there are four pairs of pants in the complement.  We will call the one or two pairs of pants that intersect $T \times \{0\}$ the \textit{bottom pants} and the one or two remaining pairs of pants the \textit{top pants}.

Consider a sequence of pants decompositions $P_1,\dots,P_n$ corresponding to a path in $\mathcal{P}(S)$.  The first pants decomposition $P_1$ determines a decomposition of the surface $S \times \{0\} \subset S \times [0,1]$.  Let $\ell \subset S \times \{0\}$ be the loop that is replaced in the first edge of the path.  We can embed a pants block in $S \times [0,1]$ so that the union of its bottom faces are identified with $F_\ell \times \{0\}$.

The union of $(S \setminus F_\ell) \times \{0\}$ and the top faces of this pants block will be a surface homeomorphic to $S$ with an induced pants decomposition.  Moreover, because the top and bottom loops of a pants block project to loops that intersect minimally, we can choose the embedding of the pants block so that the projection of this pants decomposition onto $S \times \{0\}$ will be isotopic to $P_2$.  We can continue in this way, embedding a pants block for each edge in the path.  If $P_1$ and $P_n$ have no loops in common then the union of these pants blocks will form a regular neighborhood of $S \times \{0\}$, and this can be isotoped so that their union is $S \times [0,1]$.

This construction produces a model of a surface cross an interval.  We would like to apply a similar layering construction to produce models of closed 3-manifolds.

\begin{Def}
A \textit{model decomposition} of a 3-manifold $M$ is a triple $(L, P, B)$ where $L$ is a link in $M$, $P$ is a collection of immersed pairs of pants with pairwise disjoint, embedded interiors and boundaries contained in $L$.  The set $B$ is a collection of embedded pants blocks with boundaries in $\bigcup P$ such that the marked loops in each pants block are sent into $L$.
\end{Def}

\begin{proof}[Proof of Theorem~\ref{mainthm2}]
To show that the 1-skeleton $L$ of a model decomposition $(L, P, B)$ is a hyperbolic, we will show that the complement is irreducible and atoroidal.  Each model block $B_i$ is a handlebody, and thus boundary compressible.  However, it has the structure of a surface cross and interval $F \times [0,1]$, so every compressing disk must intersect the annuli $\partial F \times [0,1]$.  Because $L \cap B_i$ contains an essential loop in each of these annuli, there is no compressing disk for $B_i$ disjoint from $L \cap B_i$.   Similarly, any incompressible annulus in $B_i$ is either boundary parallel or vertical in $F \times [0,1]$.  In the second case, it must be parallel to one of the annuli in $\partial F \times [0,1]$, and thus to a loop in $L \cap B_i$.

Let $S \subset (M \setminus L)$ be an embedded sphere, and assume that we have isotoped $S$ so as to minimize its intersection with the pairs or pants $P$.   If $S \cap P$ is not empty  then every loop bounds a disk in $S$, so some loop of $S \cap P$ will bound a disk in a block $B_i$.  However, this disk must be boundary parallel in $B_i$, so there is an isotopy of $S$ that eliminates this loop of intersection.  Thus $S \cap P$ is empty, $S$ is contained in a block $B$ and thus bounds a ball (since handlebodies are irreducible.) 

Given an embedded torus $T \subset (M \setminus L)$, a similar argument implies we can isotope $T$ to intersect $P$ in a collection of essential loops.  Thus its intersection with each block $B_i$ will be a collection of annuli.  As noted above, each annulus must be parallel to a loop in $L \cap B_i$.  As we follow the annuli around $T$, we see that they must all be parallel to the same component of $L$, and so $T$ must be boundary parallel.  Thus $M \setminus L$ contains no essential tori, and is thus hyperbolic.
\end{proof}

The construction above determines a model decomposition for $S \times [0,1]$.  By choosing a model so that the pants decompositions for $S \times \{0\}$ and $S \times \{1\}$ are related by a given homeomorphism $\phi : S \rightarrow S$, one can construct a model structure for the surface bundle with monodromy $\phi$.  By a similar construction for surfaces with boundary, one can construct model decompositions for any closed 3-manifold from an open decomposition.  However, in order to get the control needed to prove Theorem~\ref{mainthm}, we will construct a model decomposition from a generalized Heegaard splitting.

\section{Generalized Heegaard splittings}
\label{genhsect}

\begin{Def}
A handlebody is a 3-manifold homemorphic to the regular neighborhood of a graph embedded in $S^3$.  A \textit{strict compression body} $H$ is a connected 3-manifold with boundary that results from attaching 1-handles to the $S \times \{1\}$ boundary component of $S \times [0,1]$ for a compact, closed (not necessarily connected) surface $S$.  The image in $H$ of $S \times \{0\}$ is the \textit{negative boundary} $\partial_- H$ and the complement $\partial H \setminus \partial_- H$ is the \textit{positive boundary}, $\partial_+ H$.  A \textit{compression body} is a disjoint union of zero or more handlebodies and zero or more strict compression bodies.
\end{Def}

Note that we use a very general definition of compression body here.  In particular, a compression body may not be connected.  

\begin{Def}
A \textit{generalized Heegaard splitting} for a closed 3-manifold $M$ is a collection of compression bodies $H^-_1, H^+_1,\dots,H^-_k, H^+_k \subset M$ with pairwise disjoint interiors whose union is $M$ and whose boundary components coincide as follows:  For each $i$, $\partial_+ H^-_i = \partial_+ H^+_i$ and $\partial_- H^+_i = \partial_- H^-_{i+1}$.
\end{Def}

For each $i$, we will define the submanifold $M_i = H^-_i \cup H^+_i \subset M$.  The pair of handlebodies $H^-_i, H^+_i$ determines a Heegaard splitting for $M_i$ along the Heegaard surface $\Sigma_i = \partial_+ H^-_i = \partial_+ H^+_i$.  Each $\Sigma_i$ is called a \textit{thick surface} of the generalized Heegaard splitting.  The surfaces that come from the negative boundary components of the compression bodies are called \textit{thin surfaces}.

If for some $i$ there is an essential, properly embedded disk in $H^-_i$ whose boundary is disjoint from an essential disk in $H^+_i$ then following Scharlemann-Thompson~\cite{schtom}, we can replace $H^-_i, H^+_i$ with a generalized Heegaard splitting for $M_i$ whose intermediate surfaces have lower genus.  This process is called a \textit{weak reduction}.

If it is not possible to carry out a weak reduction then we will have the condition that for each $i$, every essential disk in $H^-_i$ must intersect every essential disk in $H^+_i$ in a non-empty set in $\Sigma_i$.  A generalized Heegaard splitting that satisfies this property will be called \textit{strongly irreducible}.  The main results of Scharlemann and Thompson's paper can be summarized as follows:

\begin{Thm}[Scharlemann-Thompson~\cite{schtom}]
\label{scharthomthm}
Every 3-manifold admits a strongly irreducible generalized Heegaard splitting and the thin surfaces of any strongly irreducible generalized Heegaard splitting are incompressible.
\end{Thm}

We will also need the following property of strongly irreducible generalized Heegaard splittings.  (Scharlemann's proof is for Heegaard splittings rather than generalized splittings, but the result can be generalized directly due to the fact that the thin surfaces in a strongly irreducible generalized Heegaard splitting are incompressible.)

\begin{Lem}[Scharlemann's no-nesting Lemma~\cite{schar:local}]
\label{nonestlem}
If $\ell$ is a simple closed curve in a thick surface $\Sigma_i$ of a strongly irreducible generalized Heegaard splitting such that $\ell$ bounds a disk in $M$ then $\ell$ bounds a disk contained in $H^-_i$ or $H^+_i$.
\end{Lem}

Note that if $M$ is hyperbolic then every compression body in a strongly irreducible generalized Heegaard splitting of $M$ will have no boundary component of genus one.

\section{Compression bodies and fat spines}
\label{compsect}

We will construct a model decomposition for a 3-manifold $M$ by choosing a strongly irreducible generalized Heegaard splitting for $M$ and constructing a model decomposition for each compression body in the splitting so that they match up along their boundary surfaces.  Our construction will proceed by induction, starting with models for the simplest types of compression bodies.

\begin{Def}
Given a compact, closed, orientable surface $S$, a \textit{fat spine} for $S \times [0,1]$ is a pants decomposition for $S \times \{\frac{1}{2}\}$.  A \textit{fat spine} for a genus-two handlebody is an embedded pair of pants with a pants decomposition consisting of its boundary loop and a loop in its interior, whose complement in $H$ is homeomorphic to $\partial H \times [0,1)$.
\end{Def}

Each type of fat spine is a union of pairs-of-pants whose boundary loops form a link in $H$.  Moreover, the complement of each spine is a regular neighrbohood of the boundary of $H$.  Consider a pants block $B$ that intersects a fat spine $K$ in its bottom pants and has interior disjoint from $K$.  Then the complement of $K \cup B$ will also be a regular neighborhood of $\partial H$.  After attaching this first pants block, we can attach a second pants block along its bottom pairs of pants and so on for any number of pants blocks.

\begin{Def}
A model complex in the interior of $H$ constructed in this way will be called a \textit{layered model} for $H$.
\end{Def}

Every compression body with no tori in its boundary can be constructed by attaching 1-handles to one or more of these initial types of compression bodies.  More precisely, given any compression body $H$ that is not one of these two initial types, there is a properly embedded, essential disk $D \subset H$ such that the complement of $H$ is one or two compression bodies whose positive boundary or boundaries are atoroidal and have lower genera than that of $H$.

\begin{Def}
A \textit{fat spine} for $H$ is the union of a model for each component of $H \setminus D$ and a pair of pants $F$ with one boundary loop in the boundary of each model such that the complement of the fat spine is homeomorphic to $\partial H \times [0,1)$.
\end{Def}

A fat spine can be constructed by induction for any compression body with atoroidal boundary.  The pair of pants $F$ has either one or two boundary components in the models for the smaller compression bodies, depending on whether $D$ is separating or non-separating.  % The pair of pants $F$ will be called the \textit{flat pants} in the spine.

\section{Models for compression bodies}
\label{modelcompsect}

An \textit{interior loop} of a fat spine $K$ is a loop $\ell$ in the 1-skeleton of $K$ such that some regular neighborhood of $\ell$ is contained in the pants blocks of $K$.  Any non-interior loop in the 1-skeleton is called a \textit{boundary loop}.  As with the fat spines for our initial types of compression bodies, we can layer pants blocks onto any fat spine, and the resulting model decomposition will be called a \textit{model} for $H$.  Each model for $H$ induces a pants decomposition for $\partial H$ as follows:

\begin{Def}
A pants decomposition for $\partial H$ is \textit{induced by} a model $K$ for $H$ if for each loop $\ell$ in the pants decomposition, there is an annulus $A_\ell \subset H$ with one boundary component coinciding with $\ell$, the other component coinciding with a boundary loop in $K$, and interior disjoint from $K$.  The annulus $A_\ell$ will be called a \textit{vertical annulus}.  A pants decomposition of $\partial H$ is \textit{spinal} if it is induced by a fat spine for $H$.
\end{Def}

In many cases, the union of the pants blocks in the model for $H$ will be isotopic onto $H$, and the induced pants decomposition comes from the image after this isotopy.  However, in some models (such as fat spines) the union of the pants blocks and any loose pair of pants, will not be a 3-manifold.

\begin{Lem}
Each model $K$ for $H$ induces a unique pants decomposition for $\partial H$.  
\end{Lem}

\begin{proof}
Let $\ell$ be a loop in $K$ that is the boundary of one or more pairs of pants.  A regular neighrborhood of $\ell$ is a solid torus $T$ such that $K \cap T$ consists of a collection of annuli, each with one boundary component on $\ell$ and the other in $\partial T$.  Let $A''$ be a collection of annuli with interior in the complement $T \setminus K$ such that each annulus in $A''$ has a boundary component on $\ell$ and a boundary loop in $\partial T$, and such that each component of $T \setminus K$ contains exactly one annulus of $A''$.  Construct such a collection for each loop of $K$ and let $A'$ be their union.

The complement $H \setminus K$ is homeomorphic to $\partial H \times [0,1)$ and the image of $A'$ in this product is a collection of annuli, whose closures in $\partial H \times [0,1]$ have essential boundary components in $\partial H \times \{1\}$.  Let $A \subset \partial H \times [0,1]$ be the result of extending each annulus of $A'$ to an annulus with one boundary component in $\partial H \times \{0\}$ and one in $\partial H \times \{1\}$.  By construction, the intersection of $A$ with $\partial H \times \{1\}$ is a pants decomposition.  Because the annuli are all essential, they determine the same isotopy classes of loops in $\partial H \times \{0\}$, so the intersection of $A$ with this surface is a pants decomposition.  Moreover, the images in $H$ of the annuli in $A$ show that this pants decomposition is induced by $K$.

Any annulus with a boundary component on a loop of $K$ and the other in $\partial H$ is isotopic into an annulus of $A$ within a regular neighborhood of the loop.  After this isotopy the two loops will determine the same loop in $\partial H \times \{1\}$, so the two annuli will be isotopic in $\partial H \times [0,1]$ and thus in $H$.  Thus every spine of $H$ determines a unique spinal pants decomposition.
\end{proof}

By definition, every model for a compression body $H$ contains a fat spine $K$ for $H$.  The \textit{layer number} of the model is the smallest number of model blocks not contained in $K$, over all models for $H$ and fat spines $K$ in those models.  

\begin{Lem}
\label{minlayerlem}
Let $H$ be a compression body and let $P$ be a pants decomposition for $\partial_+ H$.  If $K$ is a model for $H$ that induces $P$ on $\partial_+ H$, and has the minimal layer number among all such models then any loop in $K$ that bounds a disk in $H$ is in the interior of the fat spine.
\end{Lem}

\begin{proof}
Let $K$ be a fat spine for $H$ contained in the model.  The model blocks not in $K$ determine a path in the pants complex for $\partial_+ H$ from a spinal decomposition to $P$.  Let $P_1,\dots,P_n$ be the vertices of this path, where $P = P_1$ and $P_n$ is induced by $K$.  Let $i$ be the index such that $P_i$ is the first pants decomposition containing a loop bounding a disk in $H$.  Let $\ell \subset P_i$ be the loop bounding a disk.

By assumption, $P$ does not contain such a loop so $i > 1$.  %Moreover, the final edge in does not create a loop bounding a disk, so $i < n$.  
We will show that $P_{i-1}$ is spinal, implying that the model does not have minimal layer number.

Note that if $F_\ell$ is a punctured torus then $\ell$ is non-separating in $F_\ell$.  Moreover, if one takes two parallel copies of a disk bounded by $\ell$ and attaches them by a band in the complement of the annulus between the boundaries of two disks, the resulting disk has boundary parallel to the puncture in the torus.  Thus the loop in the boundary of $F_\ell$ already bounds a disk in $H$.  Since we chose $i$ to be the earliest edge creating such a loop, we can assume that $F_\ell$ is a four-punctured sphere.

Let $D$ be a disk with $\partial D = \ell$ and let $H' \subset H$ be the complement in $H$ of an open regular neighborhood $N$ of $D$.  Then $H'$ is either a compression body or a disjoint union of two compression bodies (depending on whether or not $D$ is separating).  

The intersection of the closure of $N$ with $H'$ is a pair of disks $E_1, E_2 \subset \partial H'$.  The restriction of the pants decomposition $P_i$ to $\partial H'$ forms a pants decomposition for the punctured surface $\partial H' \setminus (E_1 \cup E_2)$.  The two pairs of pants containing the boundary loops of this surface are the restrictions of $F_\ell$.  Each of these pairs of pants has one boundary component that bounds a disk in $\partial H'$, and the union of this pair of pants and the disk is an annulus. Let $A_1$, $A_2$ be these two annuli in $\partial H'$.  

If $A_1$ is contained in a solid torus component of $H'$ then we will define $\ell_1$ to be a core of the solid torus.  Otherwise, we will crush $A_1$ down to a single loop $\ell_1$ in $\partial H'$.  Define $\ell_2$ similarly.  Crushing $A_1$ and $A_2$ to loops (if they are not in solid tori components) induces a pants decomposition $P'$ on the boundary of the non-solid-torus component(s) of $H'$.

For each component of $H'$ that is not a solid torus, choose a model that induces the pants decomposition $P'$ on $\partial H'$.  Consider the result of attaching a pair of pants to this model for $H'$, sending two of the boundary components to the loops $\ell_1, \ell_2$.  If $\ell_1$ or $\ell_2$ is in a solid torus then the pair of pants will be attached to $H'$ along only one boundary loop.  If both $\ell_1$ and $\ell_2$ are in solid tori then $H$ is a genus-two handlebody and the resulting fat spine consists of a single pair of pants.

The induced pants decomposition of $\partial H$ is identical to $P_i$ outside of $F_\ell$.  Within $F_\ell$, the induced pants decomposition consists of a loop that intersects $\ell$ in two points.  By twisting around the disk $D$, we can make this loop into any essential loop in $F_\ell$ that intersects $\ell$ in two points.  In particular, we can make it into the loop that preceded $\ell$ in $P_{i-1}$.  Thus $P_{i-1}$ is spinal and the model for $H$ does not have the minimal layer number.
\end{proof}

\begin{Coro}
\label{minlayercoro}
Let $H$ be a compression body and let $P$ be a pants decomposition for $\partial_+ H$ such that no loop in $P$ bounds a disk in $H$.  Then there is a model $K$ for $H$ that induces $P$ on $\partial_+ H$, has minimal layer number among all such models and such that no loop in $K$ bounds a disk in $H$.
\end{Coro}

\begin{proof}
We will prove this by induction on the genus of $H$, starting with the base case of genus two.  In this case a fat spine is a single pair of pants embedded in $H$ and it is straightforward to check that no loop bounds a disk.  Moreover, if we extend this to a model that minimizes the layer number then Lemma~\ref{minlayerlem} implies that none of the new loops will bound disks in $H$.

For the induction step, assume that the Corollary holds for every genus strictly less than $g$.  Choose a model for $H$ that minimizes the layer number.  Let $K$ be the fat spine and let $H' \subset H$ be the handlebody contained in $K$ (i.e. we get the fat spine for $H$ by attaching a pair of pants to $H'$.)  The model for $H$ determines a pants decomposition for $\partial H'$ such that no loop in this pants decomposition bounds a disk in $H$.  Therefore, no loop bounds a disk in $H'$.  By the induction assumption, this implies that we can choose a model for $H'$ such that no loop bounds a disk in $H'$.  Since this model for $H'$ agrees with our original model for $H'$, we can assume that the restriction of $K$ to $H'$ is such a model.

The boundary of $H'$ is incompressible into $H \setminus H'$, so if a loop in $H'$ bounds a disk in $H$ then it bounds a disk in $H'$.  By the contrapositive, no loop in the model can bound a disk in $H$ because it does not bound a disk in $H'$.  Thus the model we've constructed contains no loop  bounding a disk in $H$.
\end{proof}

\section{Constructing the  model}
\label{constrsect}

\begin{proof}[Proof of Theorem~\ref{mainthm}]
By Scharlemann-Thompson's Theorem (Theorem~\ref{scharthomthm}), every hyperbolic 3-manifold $M$ admits a strongly irreducible generalized Heegaard splitting $\{H^-_1,H^+_1,\dots,H^-_k,H^+_k\}$ whose thin surfaces are incompressible.  Because $M$ is closed and hyperbolic, each thin surface must have genus at least two.  Because $M$ is not a lens space or a connect sum containing a lens space, each thick surface must also have genus at least two.

For each $i$, choose a fat spine $K^+$ for $H^+$ and a model $K^-$ for $H^-$ so that the pants decomposition of $\Sigma_i = \partial_+ H^-_i$ is isotopic to the decomposition induced by $K^-$ on $\Sigma_i = \partial_+ H^+_i$.  Moreover, choose $K^-_i$, $K^+_i$ so that no loop in the fat spines of $K^+_i$, $K^-_i$ bounds a disk in $H^-_i$ $H^+_i$, respectively and so that $K^-_i$ has minimal layer number over all such pairs of models.

Because $K^+_i$ is a fat spine for $H^+_i$, every boundary loop in $K^+_i$ is either disjoint from a compressing disk for $K^+_i$ or is the loop added with the final pair of pants.  Because $\Sigma_i$ is strongly irreducible, if a boundary loop of $K^+$ bounds a disk in $H^-_i$ then it must be this final loop $\ell$.  If this is the case then there is a meridian disk $D$ for $H^+_i$ that intersects the final pair of pants in an arc and intersects $\ell$ in two points.  Applying a half-twist along $D$ to $K^+_i$ replaces $\ell$ with a new loop $\ell'$ that intersects $\ell$ in two points.  Since $\ell$ bounds a disk in $H^-_i$, the construction in Lemma~\ref{minlayerlem} can be used to find a fat spine $K^+_i$ for $H^+_i$ in which $\ell'$ is the loop added with the final pair of pants.  In this case the layer number is zero, so this pair $K^-_i$, $K^+_i$ will minimize the layer number.

Thus we can assume no loop in $K^+_i$ bounds a disk in $H^-_i$.  In particular, no loop in the induced pants decomposition of $\Sigma_i$ bounds a disk in $H^-_i$, so by Lemma~\ref{minlayercoro}, we can choose $K^-_i$ so that no loop bounds a disk in $H^-_i$.  The last block layered into $H^-_+$ has top pants that project to a pair of pants in the projection of $K^+_i$.  If we remove this pants block from $K^-_i$ and layer it onto $K^+_i$, the projections of the resulting models for $H^-_i, H^+_i$ will still coincide (up to isotopy) in $\Sigma_i$.  We can transfer each subsequent model block outside the fat spine in $H^-_i$ in a similar fashion until we have a fat spine for $H^-_i$ and a model for $H^+_i$.  If this model does not have minimal layer number then we could choose a model with lower layer number, then transfer the blocks back into $H^-_i$, to reduce its layer number.  Thus this model for $H^+_i$ must have minimal layer number.  Thus by Lemma~\ref{minlayerlem}, no loop added by the layering (i.e. outside the fat spine) will bound a disk in $H^+$.  

Any loop in the interior of the fat spine of $H^-_i$ or $H^+_i$ is disjoint from a compressing disk for the compression body.  Thus if such a loop bounds a disk in the opposite compression body, the induced Heegaard splitting of $M_i$ will be weakly reducible.  Thus no loop in $K^-_i$ or $K^+_i$ bounds a disk in either compression body.  Moreover, by Scharlemann's no-nesting Lemma (Lemma~\ref{nonestlem}) these loops cannot bound disks in $M_i$.

We would like to turn the models for the compression bodies in the generalized Heegaard splitting for $M$ into a single model for $M$.  For each thin surface $\partial_- H^+_i = \partial_- H^-_{i+1}$, choose a path in the pants complex between the projections of the models we have chosen for the two compression bodies.  If there is a loop that appears in both compression bodies, choose a non-minimal path that passes through a pants decomposition that does not contain any such loop.  This path will induce a model structure for the neighborhood of $\partial_- H^+_i = \partial_- H^-_{i+1}$ between the two compression body models that matches up with the models along the intersection.  Because each thin surface is incompressible, these loops will not bound disks in $M$.

To extend the models along the thick surfaces, recall that on each $\Sigma_i$ there is a pants decomposition such that each loop is connected to a loop in $H^-_i$ and $H^+_i$ by an embedded annulus with interior disjoint from $K^\pm_i$ and its second boundary loop in $K^\pm_i$.  Let $A$ be the union of these annuli.  If $A$ contains an embedded torus $T \subset M_i$ then this embedded torus must be compressible (since $M$ is atoroidal).  If $T$ bounds a solid torus then we will throw away all the model cells inside the solid torus and crush the solid torus down to a loop. 

Otherwise, assume for contradiction that $T$ does not bound a solid torus. Let $D$ be a compressing disk for $T$ and let $X \subset M_i$ be the component of $M_i \setminus T$ that does not contain $D$.  Let $\Sigma \subset M_i$ be a regular neighborhood of the fat spine for $H^+_i$.  Let $S$ be the result of isotoping $T$ in the complement of $X$ until there is a compressing disk for $T$ disjoint from $\Sigma$ and compressing.  Then the complement of $S$ will be a ball $B$ containing $X$.  Because $\Sigma$ is strongly irreducible and $S$ is incompressible in the complement of $\Sigma$, Scharlemann's characterization of how a strongly irreducible Heegaard surface intersects a ball~\cite{schar:local} implies that $\Sigma \cap B$ must be planar and unknotted.  In particular, any loop in $\Sigma \cap B$ bounds an embedded disk in the complement of $\Sigma$.  However this is impossible because $\Sigma \cap B$ contains one or more loops in a spinal pants decomposition for $\Sigma$, which by construction do not bound a disks in $H^+_i$ or $H^-_i$. 

Thus we can crush each torus in $A$ down to a loop.  The resulting set $A'$ is homeomorphic to the product of $S^1$ with a disjoint union of trees.  Thus we can shrink each component of $A'$ down to a loop.  The annuli in $A'$ bound regions of $M_i \setminus (K^+_i \cup K^-_i)$ with the structure of a pair of pants cross an interval.  When we shrink the annuli to loops, we can isotope $K^-_i$ and $K^+_i$ into these regions, to shrink each one down to a pair of pants.  After this isotopy, the union of the models for the compression bodies is a model decomposition for $M$ such that no loop bounds a disk in $M$.
\end{proof}

\bibliographystyle{amsplain}
\bibliography{models}

%\begin{figure}[htb]
%  \begin{center}
%  \includegraphics[width=2in]{locknot.eps}
%  \caption{A graph which is knotted in a ball because of a trivial cone.}
%  \label{locknot}
%  \end{center}
%\end{figure}

\end{document}